\documentclass[11pt,leqno]{smfart}

\usepackage[latin1]{inputenc}
\usepackage[frenchb]{babel}
\usepackage{amsmath} 
\usepackage{amssymb}
\usepackage{a4wide}
\usepackage{epsfig}
\usepackage{enumerate}
\usepackage{stmaryrd}

\newtheorem{theoreme}{Th\'eor\`eme}

\renewcommand{\theequation}{\arabic{equation}}
\newcounter{hypo}
\newcounter{hyphyp}

\def\N{{\mathbb N}} 
\def\R{{\mathbb R}}

\def\oph{\mathop{\rm Op}_{h}\nolimits}

\def\Hess{\mathop{\rm Hess}\nolimits}
\newcommand{\tr}{\operatorname{tr}}
\newcommand{\one}{\operatorname{\mathbf 1}}
\def\<{\langle}
\def\>{\rangle}

\renewcommand{\theenumi}{\sl{\roman{enumi}}}
\renewcommand{\labelenumi}{\theenumi)}

\makeatletter
 \@addtoreset{equation}{section}
 \makeatother
 \renewcommand{\theequation}{\thesection.\arabic{equation}}

\begin{document}

\author[J.-F. Bony]{Jean-Fran\c{c}ois Bony}
\address[J.-F. Bony]{IMB (UMR CNRS 5251), Universit\'e Bordeaux 1, 33405 Talence, France}
\email[J.-F. Bony]{bony@math.u-bordeaux1.fr}

\author[N. Burq]{Nicolas Burq}
\address[N. Burq]{LMO (UMR CNRS 8628), Universit\'e Paris Sud 11, 91405 Orsay, France}
\email{nicolas.burq@math.u-psud.fr}

\author[T. Ramond]{Thierry Ramond}
\address[T. Ramond]{LMO (UMR CNRS 8628), Universit\'e Paris Sud 11, 91405 Orsay, France}
\email{thierry.ramond@math.u-psud.fr}

\keywords{Estimation de la résolvante, Opérateur de Schr\"{o}dinger semiclassique, trajectoires captées}

\subjclass{35J10, 35P05, 47A10, 81Q20}

\title{Minoration de la résolvante dans le cas captif}

\date{}

\begin{abstract}
Dans cette note, on démontre une minoration universelle optimale sur la norme de la résolvante tronquée pour les opérateurs de Schr\"{o}dinger semiclassiques près d'une énergie captive. En particulier, ce résultat implique que des majorations connues pour des captures hyperboliques sont optimales. La preuve repose sur un argument de X. P. Wang et la propagation en temps d'Ehrenfest des états cohérents.
\end{abstract}

\begin{altabstract}
{\bf Lower bound on the resolvent for trapped situations}
In this note, we prove an optimal universal lower bound on the truncated resolvent for semiclassical Schr\"{o}dinger operators near a trapping energy. In particular, this shows that known upper bounds for hyperbolic trapping are optimal. The proof rely on an idea of X. P. Wang, and on propagation of coherent states for Ehrenfest times.
\end{altabstract}

\maketitle

\section{Introduction}

On considère sur $\R^{n}$, $n \geq 1$, un opérateur de Schr\"{o}dinger semiclassique
\begin{equation*}
P = - h^{2} \Delta + V(x) ,
\end{equation*}
où le potentiel $V \in C^{\infty} ( \R^{n} ; \R )$ est à longue portée:
\begin{equation*}
\vert \partial_{x}^{\alpha} V (x) \vert \lesssim \< x \>^{- \sigma - \vert \alpha \vert} ,
\end{equation*}
pour un certain $\sigma >0$ et tout $\alpha \in \N^{n}$.

Le principe d'absorption limite dit alors que, pour tout $E_{0} > 0$,
\begin{equation*}
(P- (E_{0} \pm i 0))^{-1} := \lim_{\delta \searrow 0} (P - ( E_{0} \pm i \delta))^{-1} ,
\end{equation*}
existe en tant qu'opérateur de $L^{2} ( \< x \>^{2 s} d x)$ dans $L^{2} ( \< x \>^{- 2 s} d x)$, pour tout $s > 1/2$. De plus, $E \mapsto \< x \>^{-s} ( P-E)^{-1} \< x \>^{-s}$ est continu dans $]0 , + \infty [ \pm i [ 0 , + \infty [$.

On note $p (x, \xi ) = \xi^{2} + V ( x )$ le symbole semiclassique de l'opérateur $P$ et
\begin{equation*}
H_{p} = \partial_\xi p \cdot \partial_x - \partial_x p \cdot \partial_\xi = 2 \xi \cdot \partial_x - \nabla V (x) \cdot \partial_\xi ,
\end{equation*}
son champ hamiltonien. L'ensemble des trajectoires captées à l'énergie $E_{0} >0$ est défini par
\begin{equation*}
K (E_{0}) = \{ \rho \in p^{-1} (E_{0}) ; \ t \mapsto \exp (t H_{p} ) ( \rho ) \text{ reste borné en temps} \} .
\end{equation*}
C'est un ensemble compact de $T^{*} \R^{n}$. On dit que $E_{0}$ est captive si $K( E_{0} ) \neq \emptyset$ et non-captive dans le cas contraire. On notera aussi $\pi_{x} : (x , \xi ) \mapsto x$ la projection spatiale. Dans la cas non-captif, on rappelle le résultat classique suivant.

\begin{theoreme}\sl \label{a3}
Pour $E_{0} > 0$, les conditions suivantes sont équivalentes:
\begin{enumerate}[i)]
\item L'énergie $E_{0}$ est non-captive.

\item Il existe $\varepsilon >0$ tel que, pour tout $s > 1/2$,
\begin{equation} \label{a9}
\sup_{z \in [E_{0} - \varepsilon , E_{0} + \varepsilon ]} \big\Vert \< x \>^{- s} ( P -z)^{-1} \< x \>^{- s} \big\Vert \lesssim h^{-1} .
\end{equation}
\end{enumerate}
\end{theoreme}

On devrait écrire $(P - (z \pm i 0 ))^{-1}$ pour préciser qu'on prend la limite par dessus ou par dessous. Mais comme $\Vert \< x \>^{- s} ( P-z)^{-1} \< x \>^{- s} \Vert = \Vert \< x \>^{- s} ( P- \overline{z} )^{-1} \< x \>^{- s} \Vert$, ces deux limites ont même norme et cette précision est superflue.

L'implication {\sl i)} $\Rightarrow$ {\sl ii)} a été démontrée par Robert et Tamura \cite{RoTa87_01} grâce aux constructions d'Isozaki et Kitada. L'estimation \eqref{a9} a également été prouvée par Gérard et Martinez \cite{GeMa88_01} à la l'aide de la théorie de Mourre, par le deuxième auteur \cite{Bu02_01} pour des perturbations $C^{2}$ à support compact en utilisant la propagation des mesures semiclassiques et par Castella et Jecko \cite{CaJe06_01} pour des potentiels $C^{2}$ à courte portée en dimension $n \geq 3$. L'implication {\sl ii)} $\Rightarrow$ {\sl i)} du théorème \ref{a3} est due à Wang \cite{Wa87_01,Wa91_01}.

Dans le cas captif, nous démontrons ici le

\begin{theoreme}\sl \label{a1}
Pour $E_{0} > 0$, les conditions suivantes sont équivalentes:
\begin{enumerate}[i)]
\item L'énergie $E_{0}$ est captive.

\item Pour tous $\varepsilon >0$ et $\chi \in C^{\infty}_{0} ( \R^{n} )$ tel que $\chi =1$ près de $\pi_{x} K ( E_{0} )$,
\begin{equation} \label{a5}
\sup_{z \in [E_{0} - \varepsilon , E_{0} + \varepsilon ]} \big\Vert \chi ( P -z)^{-1} \chi \big\Vert \gtrsim h^{-1} \vert \ln h \vert .
\end{equation}
\end{enumerate}
\end{theoreme}

Ainsi, il y a toujours un écart d'au moins $\vert \ln h \vert$ entre les situations non-captives et captives. Ce théorème peut être démontré dans d'autres cadres. Par exemple, on peut considérer des opérateurs différentiels du second ordre, des géométries non euclidiennes voire des situations moins régulières. En fait, il suffit que la propagation des états cohérents en temps d'Ehrenfest soit valable. Du coup, on peut donner une minoration explicite de la constante qui est cachée dans l'estimation \eqref{a5} quand $h \to 0$. C'est ce qui est fait à la fin de la section \ref{a11}.

Par ailleurs, le deuxième auteur \cite{Bu04_01} à l'extérieur d'obstacles convexes (voir aussi Ikawa \cite{Ik88_01}), Christianson \cite{Ch07_01} pour une trajectoire hyperbolique, Alexandrova avec deux des auteurs \cite{AlBoRA08_01} pour un point fixe hyperbolique et Nonnenmacher et Zworski \cite{NoZw09_01} pour des captures hyperboliques à faible pression topologique ont démontré l'estimation
\begin{equation} \label{a10}
\sup_{z \in [E_{0} - \varepsilon , E_{0} + \varepsilon ]} \big\Vert \chi ( P -z)^{-1} \chi \big\Vert \lesssim h^{-1} \vert \ln h \vert .
\end{equation}
La minoration \eqref{a5} ainsi que les majorations \eqref{a10} obtenues dans les articles cités ci-dessus sont donc optimales. Ceci est également en accord avec l'intuition que les captures hyperboliques sont les captures les plus faibles. Par contre, la norme de la résolvante peut être beaucoup plus grande pour d'autres captures. C'est par exemple le cas du puits dans l'isle où
\begin{equation} \label{a2}
\sup_{z \in [E_{0} - \varepsilon , E_{0} + \varepsilon ]} \big\Vert \chi ( P -z)^{-1} \chi \big\Vert \gtrsim e^{\nu /h} .
\end{equation}
avec $\nu >0$. On sait aussi que ce type de croissance est le pire possible (voir par exemple \cite{Bu02_02}, \cite{CaVo02_01}).

Dans le théorème \ref{a1}, il suffit que $\chi$ soit strictement positif sur la projection spatiale d'une seule trajectoire semi-captée en temps positif (ou en temps négatif) de la surface d'énergie $E_{0}$. Par contre, il ne suffit pas de supposer que $\chi$ n'est pas identiquement nul. En effet, l'article de Cardoso et Vodev \cite{CaVo02_01}, transposé à notre cadre semiclassique, montre qu'il existe $R , \varepsilon >0$ tels que, pour tout $\chi \in C^{\infty}_{0} ( \R^{n} \setminus B (0,R) )$,
\begin{equation*}
\sup_{z \in [E_{0} - \varepsilon , E_{0} + \varepsilon ]} \big\Vert \chi ( P -z)^{-1} \chi \big\Vert \lesssim h^{-1}.
\end{equation*}

D'autre part, il est possible de remplacer dans \eqref{a5} la constante $\varepsilon$ par une quantité qui tend vers $0$ avec $h$. Notre preuve utilisant des fonctions qui sont localisées dans des boites de taille $h^{1/2}$, il devrait être au moins possible de prendre $\varepsilon \geq h^{\delta}$ avec $\delta < 1/2$. Par contre, l'exemple suivant montre que nécessairement $\varepsilon \geq C h$ avec $C > 0$ assez grand.

Soit $P$ un opérateur comme précédemment mais à courte portée (i.e. $\sigma > 1$) tel que, au niveau d'énergie $E_{0}$, on est dans la situation d'un puits dans l'isle ponctuel non-dégénéré (centré en $0$ pour fixer les idées). On renvoie aux hypothèses de \cite{Na89_01} pour plus de précision. La proposition 4.1 de Nakamura \cite{Na89_01} donne
\begin{equation*}
\sup_{z \in [E_{0} - \mu h , E_{0} + \mu h]} \big\Vert \< x \>^{- s} ( P -z)^{-1} \< x \>^{- s} \big\Vert \lesssim h^{-1} .
\end{equation*}
pour tout $0 \leq \mu < \tr (\Hess V(0) )/2$. Cette majoration est à comparer avec \eqref{a2} qui a lieu sous les mêmes hypothèses. Par conséquent, le fait que $E_{0}$ soit une énergie captive n'implique aucune minoration particulière de la résolvante en $E_{0}$ ou même, plus généralement, dans certains voisinages de taille $h$ de $E_{0}$.

\section{Preuve du théorème \ref{a1}} \label{a11}

Il suffit de démontrer que {\sl i)} implique {\sl ii)} d'après le théorème \ref{a3}. Nous suivons la stratégie que Wang \cite{Wa91_01} a développé pour prouver l'implication {\sl ii)} $\Rightarrow$ {\sl i)} du théorème \ref{a3}. Considérons $K (h) \geq 0$ défini par
\begin{equation*}
h^{-1} K(h) = \sup_{z \in [E_{0} - \varepsilon , E_{0} + \varepsilon ]} \big\Vert \chi ( P -z)^{-1} \chi \big\Vert .
\end{equation*}
Par le principe d'absorption limite, cette quantité est finie. D'après la régularité de Kato (on renvoie à la section XIII.7 de \cite{ReSi78_01} pour une présentation de cette théorie), l'identité précédente implique
\begin{equation*}
\int_{\R} \big\Vert \chi \one_{[E_{0} - \varepsilon , E_{0} + \varepsilon ]} (P) e^{i s P} u \big\Vert^{2} d s \leq 8 h^{-1} K(h) \Vert u \Vert^{2} .
\end{equation*}
En particulier, pour $\varphi \in C^{\infty}_{0} ( [E_{0} - \varepsilon , E_{0} + \varepsilon ] ; [0,1] )$ avec $\varphi (E_{0} )=1$, on a
\begin{equation} \label{a4}
\int_{\R} \big\Vert \chi \varphi (P) e^{i t P /h} u \big\Vert^{2} d t \leq 8 K(h) \Vert u \Vert^{2} .
\end{equation}

Par le calcul fonctionnel des opérateurs pseudodifférentiels (voir le chapitre 8 du livre de Dimassi et Sj\"{o}strand \cite{DiSj99_01}), $\varphi (P) \chi^{2} \varphi (P)$ est un opérateur pseudodifférentiel semiclassique de symbole $a (x, \xi ;h) \in S_{h} (1)$ (la classe des fonctions bornées ainsi que leur dérivées uniformément par rapport à $h$). De plus,
\begin{equation*}
a (x, \xi ;h) = a_{0} (x, \xi ) + h S_{h} (1) \qquad \text{ avec } \qquad a_{0} (x, \xi ) = \chi^{2} (x) \varphi( p ( x , \xi ))^{2} .
\end{equation*}
Ainsi, \eqref{a4} s'écrit
\begin{equation} \label{a6}
\int_{\R} \big( \oph (a) e^{i t P /h} u , e^{i t P /h} u \big) d t \leq 8 K(h) \Vert u \Vert^{2} .
\end{equation}

On utilise maintenant la propagation des états cohérents en temps d'Ehrenfest. Soit $u_{\rho}$ un état cohérent gaussien normé centré autour de $\rho \in T^{*} \R^{n}$. D'après la proposition 5.1 de Bouzouina et Robert \cite{BoRo02_01},
\begin{equation*}
\lim_{h \to 0} \big( \oph ( a ) e^{i t P /h} u_{\rho} , e^{i t P /h} u_{\rho} \big) = a_{0} ( \exp ( t H_{p} ) ( \rho ) ) ,
\end{equation*}
uniformément pour $\vert t \vert \leq \vert \ln h \vert (1 - \varepsilon ) / (2 \Gamma)$ où $\Gamma$ est une constante explicite liée à stabilité du système classique et $\varepsilon >0$ est quelconque. On fixe alors $\rho \in K (E_{0})$ qui n'est pas vide par hypothèse. Ainsi $\exp ( t H_{p} ) ( \rho )$ reste en tout temps dans $K (E_{0})$. Et comme $a_{0} =1$ sur cet ensemble, pour tout $h$ assez petit,
\begin{equation} \label{a8}
\big( \oph ( a ) e^{i t P /h} u_{\rho} , e^{i t P /h} u_{\rho} \big) \geq (1- \varepsilon ) ,
\end{equation}
uniformément pour $\vert t \vert \leq \vert \ln h \vert (1 - \varepsilon ) / (2 \Gamma)$. En combinant \eqref{a6} et \eqref{a8}, il vient
\begin{equation*}
\vert \ln h \vert (1 - \varepsilon )^{2} / ( 8 \Gamma ) \leq K (h) ,
\end{equation*}
ce qui finit la preuve et donne une minoration explicite dans la limite $h \to 0$.

\bibliographystyle{smfplain}
\providecommand{\bysame}{\leavevmode ---\ }
\providecommand{\og}{``}
\providecommand{\fg}{''}
\providecommand{\smfandname}{et}
\providecommand{\smfedsname}{\'eds.}
\providecommand{\smfedname}{\'ed.}
\providecommand{\smfmastersthesisname}{M\'emoire}
\providecommand{\smfphdthesisname}{Th\`ese}

\end{document}